\documentclass[a4paper,11pt]{article}
\usepackage{url}
\usepackage{latexsym}
\usepackage{amssymb}
\usepackage{amsmath}
\usepackage{amsthm}
\usepackage{times}
\usepackage{xspace}
\usepackage[dvips]{graphicx}

\usepackage{geometry}
\newtheorem{remark}{Remark}
\newtheorem{claim}{Claim}
{\bfseries}{\itshape}
{\bfseries}{\itshape}
{\bfseries}{\itshape}
{\bfseries}{\itshape}
\newtheorem{theorem}{Theorem}{\bfseries}{\itshape}
\newtheorem{lemma}{Lemma}{\bfseries}{\itshape}
\newtheorem{definition}{Definition}{\bfseries}{\itshape}

\newcommand {\bbox}{\rule{0.6em}{0.6em}}
\newcommand{\boxi}{\ensuremath{\mathrm{box}}}
\newcommand{\cubi}{\ensuremath{\mathrm{cub}}}

\title{An upper bound for Cubicity in terms of Boxicity}
\author{L. Sunil Chandran
\thanks{Dept. of Computer Science and Automation, Indian Institute of Science, Bangalore--560012, India. Email : 
{\tt sunil@csa.iisc.ernet.in}}
\and K. Ashik Mathew
\thanks{Dept. of Computer Science and Automation, Indian Institute of Science, Bangalore--560012, India. Email : 
{\tt ashik@csa.iisc.ernet.in}}
}
\date{}

\begin{document}
\maketitle
\thispagestyle{empty}

\begin{abstract}
An axis-parallel \emph{b}-dimensional box is a Cartesian product
 $R_1 \times R_2 \times \ldots \times R_b$ where each $R_i$ (for $1 \leq
 i \leq b$) is a  closed interval of the form $[a_i,b_i]$ on the real line. 
The boxicity of any graph $G$, $\boxi(G)$ is the minimum positive integer $b$ such that G can be represented as the intersection graph of axis parallel $b$-dimensional boxes. 
A $b$-dimensional cube is a Cartesian product $R_1 \times R_2\times \ldots \times R_b$, where each $R_i$ \mbox{(for $1 \leq i \leq b$)} is a closed interval of the form [$a_i$,$a_i$+1] on the real line.
When the boxes are restricted to be axis-parallel cubes in \emph{b}-dimension, the minimum dimension \emph{b} required to represent the graph
 is called the cubicity of the graph (denoted by $\cubi(G)$). In this paper
 we prove that \mbox{$\cubi(G)\leq \lceil \log n \rceil  \boxi(G)$} where $n$
 is the number of vertices in the graph. We also show that this upper bound is tight. 
\end{abstract}
\section{Introduction}
Let $\mathcal{F}=\{S_x \subseteq \emph{U}: x\in V\}$ be a family of subsets of \emph{U}, where \emph{V} is an index set. The intersection graph $\Omega$($\mathcal{F}$) of $\mathcal{F}$ has $V$ as vertex set, and two distinct vertices \emph{x} and \emph{y} are adjacent if and only if $S_x \cap S_y \neq \emptyset$. Representations of graphs as the intersection graphs of various geometric objects is a well-studied area in graph theory. A prime example of a graph class defined in this way is the class of interval graphs. \\
\begin{definition}
A graph $I(V,E)$ is an interval graph if and only if there exists a function $\Pi$ which maps each vertex $u \in V$ to a closed interval of the form $[l(u),r(u)]$ on the real line such that $(u,v)\in E$ if and only if  $\Pi (u) \cap \Pi (v) \neq \emptyset $. We will call $\Pi$ an interval representation of $I(V,E)$.
\end{definition}
\begin{definition}
 An indifference graph is an interval graph which has an interval representation in which each of the intervals is of the same length. We will call such an interval representation a unit interval representation of the graph.
\end{definition}The indifference graphs are also known as unit interval graphs. See Chapter 8 of \cite{Golu} for more information on interval graphs and indifference graphs.\\

Motivated by theoretical as well as practical considerations, graph theorists have tried to generalize the concept of interval graphs in many ways. In many cases, representation of a graph as the intersection graph of a family of geometric objects, which are generalizations of intervals, is sought. Concepts such as boxicity and interval number are examples.\\

\noindent In this paper we only consider simple, finite, undirected graphs.
$V(G)$ and $E(G)$ denote the set of vertices and the set of edges of $G$, respectively.   

\begin{definition}
For a graph $G$, $\boxi(G)$ is the minimum positive integer $b$ such that $G$ can be represented as the intersection graph of axis-parallel $b$-dimensional boxes. Here a \emph{b}-dimensional box is a Cartesian product $R_1 \times R_2 \times \ldots \times R_b$ where each $R_i$ (for $1 \leq i \leq b$) is defined to be a closed interval of the form $[a_i,b_i]$ on the real line. 
The boxicity of a complete graph is defined to be $0$.
\end{definition}

\begin{definition}
The cubicity of a graph $G$, $\cubi(G)$ is the minimum positive integer $b$ 
such that $G$ can be represented as the intersection graph of axis-parallel
 $b$-dimensional cubes. Here a $b$-dimensional cube is a Cartesian product
 $R_1 \times R_2\times \ldots \times R_b$, where each $R_i$ 
\mbox{(for $1 \leq i \leq b$)} is a closed interval of the form $[a_i,a_i+1]$
 on the real line.
 The cubicity of a complete graph is defined to be $0$.
\end{definition}

The following observation is easy to make. A $1$-dimensional box is 
a closed interval on the real line and thus graphs of boxicity 1 are exactly
 the interval graphs. Similarly, the graphs with cubicity 1 are the 
indifference graphs.

\begin{lemma}
[{\bf Roberts\cite{Roberts}}]
\label {BoxicityIntervalgraphCharacterizationLemma}
 Given a graph $G$,the minumum positive integer 
$b$ such that there exist interval graphs $G_1, G_2, \ldots G_b$ 
with  $V(G)=V(G_i)$ for $1\leq i \leq b$ and  
satisfying $E(G)=E(G_1) \cap E(G_2) \cap \ldots E(G_b)$ is equal to 
$\boxi(G)$.
\end{lemma}

\begin{lemma}[{\bf Roberts\cite{Roberts}}]
\label {CubicityIndifferenceGraphCharacterizationLemma}
 Given a graph  $G$,the minumum positive integer 
$b$ such that there exist indifference graphs $G_1, G_2, \ldots G_b$ 
with $V(G)=V(G_i)$ for $1\leq i \leq b$ and 
satisfying $E(G)=E(G_1) \cap E(G_2) \cap \ldots E(G_b)$ is equal to 
$\cubi(G)$.
\end{lemma}
The concepts of cubicity and boxicity were introduced by F.S. Roberts 
\cite{Roberts}.  They find applications in niche overlap in ecology and 
in solving problems of fleet maintanence in operations research.
 (See \cite {CozRob}.) It was shown by Cozzens \cite{Coz} that computing 
the boxicity of a graph is an NP-hard problem. Later, this was improved by 
Yannakakis\cite{Yan1} , and finally by Kratochvil\cite{Kratochvil} who showed
 that deciding whether the boxicity of a graph is at most 2 itself is an
 NP-complete problem. The complexity of finding the maximum independent set
 in bounded boxicity graphs was considered by \cite {ImaiAsano, Fowler}. 
Some NP-hard problems are known to be either polynomial time solvable 
or have much better approximation ratio on low boxicity graphs. For example, 
the max-clique problem is polynomial time solvable on bounded boxicity graphs
 and the maximum independent set problem has $\log n$ approximation ratio
 for graphs with boxicity 2 \cite{Agarwal98, Berman2001}. 

There have been many attempts to find the cubicity and boxicity of graphs with
 special structures. In his pioneering work, F.S. Roberts\cite{Roberts} proved
 that the boxicity of a complete \emph{k}-partite graph (where each part has 
 at least 2 vertices)  is \emph{k}. He also proved that the cubicity of 
any graph can not be greater than $\lfloor 2n/3 \rfloor$ and the boxicity 
 cannot be greater than $\lfloor n/2 \rfloor$.
Scheinerman\cite {Scheiner} showed that the boxicity of outer planar graphs 
is at most 2. Thomassen\cite {Thoma1} proved that the boxicity of planar graphs
 is bounded above by 3. The boxicity of split graphs is investigated by 
Cozzens and Roberts\cite{CozRob}. Chandran and Sivadasan\cite{BoxCubHypTech} proved 
that the cubicity of the $d$-dimensional hypercube $H_d$ is
 $\theta (\frac {d} {\log d})$. They also proved that for 
any graph $G$, $\boxi(G) \leq tw(G)+2$ where $tw(G)$ is the treewidth of $G$ \cite{CN05}. 
This in turn throws light on the boxicity of various other graph classes. 
Roberts and Cozzens proposed a theory of dimensional properties, attempting
to generalize the concepts of cubicity and boxicity \cite {CozRob89}. These concepts
were further developed by Kratochvil and Tuza \cite {KratoTuza}.

Researchers have also tried to generalize or extend the concept of boxicity 
in various ways. The poset boxicity \cite{TroWest} , the rectangular number
 \cite{ChangWest} , grid dimension \cite{Bellantoni}, 
circular dimension\cite{Feinberg} and the boxicity of digraphs\cite{ChangWest1}
 are some examples.

\section{Our Results}
It is easy to see that for any graph $G$, $\boxi(G) \le cubi(G)$. 
In this paper we prove the following theorem:
\begin{theorem}
For a graph $G$ on $n$ vertices, 
$\cubi(G) \le  \lceil \log {n} \rceil \boxi(G)$.
Moreover, this upper bound is tight.
\end {theorem}

\subsection{Consequences of our result}

The upper bound that we developed should be useful in many cases where a 
bound for one of the two quantities (boxicity and cubicity) is already known. 
Combining our theorem with previously known upper bounds for boxicity, we
get  various upper bounds for  cubicity, which we list in the following
table.  Here $n$ denotes the number of vertices in the graph, 
$tw = treewidth(G)$ is the treewidth of $G$, 
 $ \Delta = 
\Delta(G)$ is
the maximum degree and $\omega = \omega(G)$ is the clique number, i.e. 
the number
of vertices in the biggest clique in $G$. 
Each of the references given corresponds to the paper in which the
corresponding upper bound for boxicity was proved.

~~~

\begin{tabular}{|c|c|c|}
\hline
Graph Class  & Upper bound for                & Upper bound for \\
             &     $\boxi(G)$                           &  $\cubi(G)$ \\
\hline
\hline

Chordal Graphs\cite{CN05}& $\omega +1$& $(\omega+1)\log n$\\
&$\Delta+2$ &$(\Delta+2)\log n$\\
\hline
Circular Arc Graphs\cite{CN05}& $2\omega +1 $&$(2\omega +1) \log n$\\
&$2\Delta+3$ &$(2\Delta +3)\log n$\\
\hline
AT-Free Graphs\cite{CN05} & $3\Delta $ & $(3\Delta )\log n$\\
\hline
Co-comparability graphs\cite{CN05} & $(2 \Delta +1)$ & $(2 \Delta +1)\log n$\\
\hline
Permutation Graphs\cite{CN05} & $(2\Delta+1)$ & $(2\Delta +1)\log n$\\
\hline
Planar Graphs\cite{Kratochvil} & $3$ & $3\log n$\\
\hline
Series Parallel Graphs\cite{CRB1} & $3$ & $3\log n$\\
\hline
Outer Planar Graphs\cite {Scheiner} & $2$ & $2\log n$\\
\hline
Any Graph\cite{CN05} & $tw+2$ & $(tw+2)\log n$\\
\hline
Any Graph\cite{BoxRandConf} &  $(\Delta + 2) \log n$  &  $(\Delta + 2)\log^2 n$\\
\hline
\end{tabular}
\subsubsection{Algorithmic Consequences}
Our proof provides an $O(n^2\log n)$ algorithm to represent any interval graph
 $G$ (on $n$ vertices) into a $\log n$-space as the intersection graph 
of $n$ axis parallel $\log n$-dimensional cubes, when the interval representation of $G$ is given. 
Also follows from this, a polynomial time algorithm to translate any given box representation of a graph in a $b$-dimensional space to a cube representation in $b\log n$-dimensional space.

\section {Proof of our Theorem}

\begin{lemma}[{\bf Roberts\cite{Roberts}}]
\label{IntersectionBasedUpperboundLemma}
Let $G$ be a graph and let 
$G_1, G_2, \cdots, G_j$ be graphs such that 
${\mathbf (1)} $  $V(G)=V(G_p)$ for $1 \leq p \leq j$   {\bf and}   
${\mathbf (2)}$  $E(G)=E(G_1) \cap E(G_2) \cap \ldots E(G_j)$. Then 
$\cubi(G) \leq \cubi(G_1)+\cubi(G_2)+ \ldots + \cubi(G_j)$.
\end{lemma}
\begin{lemma}
\label {IntervalGraphSufficientLemma}
Let $r(n)$ denote the largest real number such that there exists
a non-complete  graph $G$ $($i.e. a graph $G$
 such that  $\boxi(G) > 0)$ on $n$ vertices such that
 $\cubi(G) = r(n) \boxi(G)$. 
Then, there exists an interval graph $G'$ on $n$ vertices such that
$\cubi(G') = r(n)$. 
\end{lemma}
\begin{proof}
Let $G$ be a graph on $n$ vertices such that $\boxi(G)=b$ and 
$\cubi(G)= b \cdot r(n) $. Then by Lemma 
\ref {BoxicityIntervalgraphCharacterizationLemma}, there exists 
interval graphs $G_1, G_2, \cdots, G_b$ such that $V(G_i) = V(G)$ for
$1 \le i \le b$ and 
$E(G)=E(G_1) \cap E(G_2) \cap \ldots E(G_b)$. 
By Lemma \ref {IntersectionBasedUpperboundLemma}, $ r(n)\cdot b = \cubi(G) 
\le \sum_{i=1}^b \cubi(G_i)$. 
It follows that there exists at least one $i$, ($1 \le i \le b$) such that
$\cubi(G_i) \ge r(n)$. Recallin that $G_i$ is a (non-complete)
 interval graph and
thus $\boxi(G_i) = 1$ we have $\cubi(G_i) \ge r(n) \cdot \boxi(G_i)$.
From the definition of $r(n)$, it follows that $\cubi(G_i) = r(n) 
\cdot \boxi(G_i) = r(n)$, as required. 
\end{proof}

\begin {lemma}
\label {IntervalgraphOrderingLemma}
For every interval graph $G$ on $n$ vertices, there exists an ordering 
$f: V(G) \rightarrow
\{1,2,\cdots,n\}$ of its vertices such that if 
$u,v,w \in V(G)$ satisfy $f(u) < f(w) < f(v)$ and $(u,v) \in E(G)$ then 
$(u,w) \in E(G)$, also.
\end{lemma}
\begin{proof}
Consider an interval representation of $G$ and order the vertices in 
the non-decreasing order of the left
 end-points of the intervals. It is easy to verify that this order satisfies the required property.
\end{proof}

\noindent {\bf Proof of Theorem 1}

~~~~

By Lemma \ref {IntervalGraphSufficientLemma}, it is enough to show that for 
any interval 
graph $G$ on $n$ vertices, $\cubi(G) \le \lceil \log n \rceil$.
Let us first assume that $n = 2^k$ for a positive integer $k$. (We will
take care of the remaining case in the end.) Then by Lemma \ref {CubicityIndifferenceGraphCharacterizationLemma},
we only have to show that there exists $k$ indifference graphs
$I_1, I_2, \cdots, I_k$ such that $V(I_i) = V(G)$  for $1 \le i \le k$
and $E(G) =
\bigcap_{i=1}^{k} E(I_i)$. 
Let $f$ be an ordering of $V$  as described in Lemma \ref {IntervalgraphOrderingLemma}.
First we define  $k+1$ different partitions $\mathcal {P}_1, \mathcal {P}_2,
\cdots, \mathcal {P}_{k+1}$ of $V$ as follows: 
 $$\mathcal{P}_i=\{ S_1^i,S_2^i,\ldots S_{m_i}^i \}, 
\mbox { {\bf where}   $S_j^i = \{v \in V : (j-1)2^{i-1} +1 
\le f(v) \le j2^{i-1} \}$ }$$ 
The reader can easily verify that for each $i$, $1 \le i \le k+1$,
$\mathcal {P}_i$ defines a valid partition of $V$ i.e., 
$\bigcup_{j} S_j^i = V$ and $S_a^i \cap S_b^i = \emptyset$ for $a \ne b$.
 Moreover for partition
$\mathcal {P}_i$  all  blocks have same cardinality, i.e.
$|S_j^i| = 2^{i-1}$. Moreover $m_i = 2^{k-i+1}$.  For $i \le k$, $m_i$ is
an clearly an even number. The partition $\mathcal {P}_{k+1}$ contains only 
one block, namely $S_1^{k+1} = V$.

\noindent For $1 \le i \le k$, we construct the indifference graph $I_i$ 
based on the 
partition $\mathcal {P}_i$. Let 
$$A_i = S_1^i \cup S_3^i \cup \ldots \cup S_{{m_i}-1}^i
 \mbox { \  and  \ } B_i = S_2^i \cup S_4^i \cup \ldots \cup S_{m_i}^i$$
Clearly $(A_i,B_i)$ is a partition of $V$. 
Now we define the indifference graph $I_i$ by defining its unit interval 
representation $\Pi_i$ as follows:

~~~~

\noindent For $v \in B_i$: \hspace {1.5cm}
$\Pi_i (v)=[n+f(v),2n+f(v)] $.

\noindent For $v \in A_i$, if $N(v) \cap B_i = \emptyset$: \hspace {1cm}
   $\Pi_i(v) = [0,n]$.

\noindent For  $v \in A_i$, if  $N(v) \cap B_i \ne \emptyset$: (Let
     $t=max_{x \in N(v) \cap B_i} f(x)$.) \hspace {1cm}
        $\Pi_i (v)=[t,n+t]$

\begin{claim}
$E(I_i) \supseteq E(G)$ for $1\leq i \leq k$
\end{claim}

\noindent Let $(u,v) \in E(G)$. We only have to  consider the following 
 three cases.

\noindent {\bf {Case 1:}} $u\in A_i$ and $v \in A_i$. Then  $\Pi_i (u) \cap
 \Pi_i (v) \neq \emptyset$ since the point corresponding to $n$ on the real 
line is a member of both  $\Pi_i (u)$ and $\Pi_i (v)$.

\noindent {\bf {Case 2:}}$u\in B_i$ and $v \in B_i$. Here also $\Pi_i (u)
 \cap \Pi_i (v) \neq \emptyset$ since the point corresponding to $2n$ on 
the real line is a member of both  $\Pi_i (u)$ and $\Pi_i (v)$.

\noindent {\bf {Case 3:}} $u \in A_i$ and $v \in B_i$. In this case, 
let $z=\max(f(x):x \in N(u) \cap B_i)$. Now, $f(v) \leq z$, 
since $v \in N(u) \cap B_i$. Now recall that $\Pi_i (v)= [n+f(v), 2n+f(v)]$ 
and $\Pi_i (u)=[z,n+z]$. Clearly, the point corresponding to $n+z$ 
on the real line belongs to  both  $\Pi_i (u)$ and $\Pi_i (v)$, and thus $\Pi_i (u) \cap \Pi_i (v) \neq \emptyset$. $\bbox$

\begin{claim}
If $(u,v) \notin E(G)$ then there exists an $i$, $1 \leq i \leq k$ 
such that $(u,v) \notin E(I_i)$.
\end{claim}

Let $t$ be the largest integer such that for  $1 \le i \leq t$, 
$u$ and $v$ belong to different blocks of the partion $\mathcal{P}_i$,
 i.e. if $1\le i \leq t$ and  $u \in S_a^i$ and $v \in S_b^i$, 
then $a \neq b$. Clearly such a $t$ exists and in fact $t \le k$,
 since $\mathcal {P}_{k+1}$ contains only one block.
 Without loss of generality, assume that $f(u) < f(v)$.  
 We claim that if $u \in S_a^t$ and $v \in S_b^t$ then $b=a+1$,
 where $a$ is an odd number. To see this notice that
 by the definition of $t$, $u$ and $v$ belong to the same
block in $P_{t+1}$ and  if $u,v \in S_c^{t+1}$ then clearly
 $u \in S_a^t$ and $v \in S_b^t$, where $a=2(c-1)+1$ and $b=2(c-1)+2$.

 Now we will show  that $(u,v) \notin E(I_t)$. To see this, first observe that $u \in A_t$ and $v \in B_t$ since $u \in S_a^t$ where $a$ is an odd number and $v \in S_{a+1}^t$ where $a+1$ is an even number. 
If $N(u) \cap B_t = \emptyset$, clearly $(u,v) \notin E(I_t)$,
since in that case $\Pi_t(u) = [0, n]$ and $\Pi_t(v) = [n + f(v), 2n + f(v)]$
and these two intervals do not intersection. So, we can assume that
$N(u) \cap B_t \ne \emptyset$. 
Now, let $w \in B_t$
be such that  $f(w)=\max(f(x):x \in N(u) \cap B_t)$. 
We claim that $f(w) < f(v)$. Suppose  not. Then clearly $f(u) < f(v) < f(w)$.
Now by Lemma \ref {IntervalgraphOrderingLemma}, $(u,v) \in E(G)$, 
since $(u,w) \in E(G)$, contradicting the 
assumption that $(u,v) \notin E(G)$ . Now, recall that
$\Pi_t (u)=[f(w),n+f(w)]$ and $\Pi_t (v)=[n+f(v), 2n+f(v)]$. 
Since $f(w) < f(v)$  we have $\Pi_t (u) \cap \Pi_t (v) = \emptyset$ and 
thus $(u,v) \notin E(I_t)$.   $\bbox$

From Claim 1 and Claim 2 we have, 
$E(G)=E(I_1) \cap E(I_2) \cap \ldots E(I_k)$ as required. 
So by Lemma \ref {CubicityIndifferenceGraphCharacterizationLemma}, 
$\cubi(G) \leq k=\log n$.  
If $2^{k-1} < |V| <  2^k$, then add $2^k-|V|$
isolated vertices to the graph. Note that this will not change the cubicity or boxicity of the graph. Moreover  $\lceil \log n \rceil = k$, and 
the result follows.   

Finally the tightness of our result can be verified by considering the star graph on $n$ vertices, $S(n)$. (Note: The star graph $S(n)$ is the complete 
bipartite graph $K_{1,n-1}$, with a single on one side and the remaining
$n-1$ vertices on the other side.)
Its  boxicity equals $1$, since it is an interval graph.
 It is also known that \cite{Roberts} 
$\cubi(S(n))= \lceil \log (n-1) \rceil$. 
Note that when $n \neq 2^k +1$, we have $\lceil \log (n-1) \rceil = \lceil \log n \rceil$ and thus our upper bound is  tight.
$\Box$

\begin{remark}
The $k$ indifference graphs that we constructed all have a diameter less 
than or equal to $2$. Thus it follows from our proof that the edge set
of any inteval
graph can be represented as the intersection of the edge sets of
at most  $\lceil \log n \rceil$
indifference graphs of diameter at most $2$. 
\end{remark}


\begin{thebibliography}{10}

\bibitem{Agarwal98}
P.~K. Agarwal, M.~van Kreveld, and S.~Suri.
\newblock Label placement by maximum independent set in rectangles.
\newblock {\em Comput. Geom. Theory Appl.}, 11:209--218, 1998.

\bibitem{Bellantoni}
S.~Bellantoni, I.~Ben-Arroyo Hartman, T.~Przytycka, and S.~Whitesides.
\newblock Grid intersection graphs and boxicity.
\newblock {\em Discrete mathematics}, 114(1-3):41--49, April 1993.

\bibitem{Berman2001}
P.~Berman, B.~DasGupta, S.~Muthukrishnan, and S.~Ramaswami.
\newblock Efficient approximation algorithms for tiling and packing problems
  with rectangles.
\newblock {\em J. Algorithms}, 41:443--470, 2001.

\bibitem{CRB1}
Ankur Bohra, L.~Sunil Chandran, and J.~Krishnam Raju.
\newblock Boxicity of series paralel graphs.
\newblock To appear in Discrete mathematics, 2005.

\bibitem{BoxRandConf}
L.~Sunil Chandran and N.~Sivadasan.
\newblock Geometric representation of graphs in low dimension.
\newblock Accepted in the 12th Annual International Computing and Combinatorics
  Conference to be held in Taipei, Taiwan, August 2006.

\bibitem{BoxCubHypTech}
L.~Sunil Chandran and N.~Sivadasan.
\newblock On the boxicity and cubicity of hypercubes.
\newblock Submitted. Available at http://aps.arxiv.org/abs/math.CO/0605246.

\bibitem{CN05}
L.~Sunil Chandran and Naveen Sivadasan.
\newblock Treewidth and boxicity.
\newblock Submitted, Available at http://arxiv.org/abs/math.CO/0505544.

\bibitem{ChangWest}
Y.~W. Chang and Douglas~B. West.
\newblock Rectangle number for hyper cubes and complete multipartite graphs.
\newblock In {\em 29th {SE} conf. Comb., Graph Th. and Comp., Congr. Numer.
  132(1998)}, 19--28.

\bibitem{ChangWest1}
Y.~W. Chang and Douglas~B. West.
\newblock Interval number and boxicity of digraphs.
\newblock In {\em Proceedings of the 8th International Graph Theory Conf.},
  1998.

\bibitem{Coz}
M.~B. Cozzens.
\newblock Higher and multidimensional analogues of interval graphs.
\newblock Ph. D thesis, Rutgers University, New Brunswick, NJ, 1981.

\bibitem{CozRob}
M.~B. Cozzens and F.~S. Roberts.
\newblock Computing the boxicity of a graph by covering its complement by
  cointerval graphs.
\newblock {\em Discrete Applied Mathematics}, 6:217--228, 1983.

\bibitem{CozRob89}
M.~B. Cozzens and F.~S. Roberts.
\newblock On dimensional properties of graphs.
\newblock {\em Graphs and Combinatorics}, 5:29--46, 1989.

\bibitem{Feinberg}
Robert~B. Feinberg.
\newblock The circular dimension of a graph.
\newblock {\em Discrete mathematics}, 25(1):27--31, 1979.

\bibitem{Fowler}
R.~J. Fowler, M.~S. Paterson, and S.~L. Tanimoto.
\newblock Optimal packing and covering in the plane are {NP}--complete.
\newblock {\em Information Processing letters}, 12(3):133--137, 1981.

\bibitem{Golu}
Martin~C Golumbic.
\newblock {\em Algorithmic Graph Theory And Perfect Graphs}.
\newblock Academic Press, New York, 1980.

\bibitem{ImaiAsano}
H.~Imai and T.~Asano.
\newblock Finding the connected component and a maximum clique of an
  intersection graph of rectangles in the plane.
\newblock {\em Journal of algorithms}, 4:310--323, 1983.

\bibitem{Kratochvil}
J.~Kratochvil.
\newblock A special planar satisfiability problem and a consequence of its
  {NP}--completeness.
\newblock {\em Discrete Applied Mathematics}, 52:233--252, 1994.

\bibitem{KratoTuza}
J.~Kratochvil and Z.~Tuza.
\newblock Intersection dimensions of graph classes.
\newblock {\em Graphs and Combinatorics}, 10:159--168, 1994.

\bibitem{Roberts}
F.~S. Roberts.
\newblock {\em Recent Progresses in Combinatorics}, chapter On the boxicity and
  Cubicity of a graph, pages 301--310.
\newblock Academic Press, New York, 1969.

\bibitem{Scheiner}
E.~R. Scheinerman.
\newblock Intersectin classes and multiple intersection parameters.
\newblock Ph. D thesis, Princeton University, 1984.

\bibitem{Thoma1}
C.~Thomassen.
\newblock Interval representations of planar graphs.
\newblock {\em Journal of combinatorial theory, Ser {B}}, 40:9--20, 1986.

\bibitem{TroWest}
W.~T. Trotter and Jr. Douglas~B. West.
\newblock Poset boxicity of graphs.
\newblock {\em Discrete Mathematics}, 64(1):105--107, March 1987.

\bibitem{Yan1}
Mihalis Yannakakis.
\newblock The complexity of the partial order dimension problem.
\newblock {\em {SIAM} Journal on Algebraic Discrete Methods}, 3:351--358, 1982.

\end{thebibliography}

\end{document}